\documentclass[12pt, a4paper]{amsart}
\usepackage[utf8]{inputenc}
\usepackage{amsxtra,amssymb,amsthm,amsmath,amscd,mathrsfs, epsfig, eufrak}
\usepackage{amscd, amsmath, mathrsfs, amssymb, amsthm, amsxtra, bbding, epsfig, eucal, eufrak, graphicx, latexsym, mathrsfs, mathbbol, bbold}
\usepackage[all]{xy}

\usepackage{tikz}

\usepackage[normalem]{ulem}
\usepackage{soul}
\usepackage{color}

\setstcolor{red}

\makeatletter
\@namedef{subjclassname@2010}{%
\textup{2010} Mathematics Subject Classification}
\makeatother

\def \d {\,{\rm d}}

\theoremstyle{plain}

\theoremstyle{remark}

\numberwithin{equation}{section}

\addtocounter{footnote}{1}

\begin{document}


\title[\tiny Waring's problem involving D.H. Lehmer numbers]
{Waring's problem involving D.H. Lehmer numbers}
\author[\tiny Yang Qu \& Rong Ma]{Yang Qu \& Rong Ma}

\address{%
Rong Ma
\\
School of Mathematics and Statistics
\\
Northwestern Polytechnical University
\\
Xi'an
\\
Shaanxi 710072
\\
China}
\email{marong0109@163.com}

\address{%
Yang Qu
\\
School of Mathematics and Statistics
\\
Northwestern Polytechnical University
\\
Xi'an
\\
Shaanxi 710072
\\
China}
\email{1993845352@qq.com}

\date{\today}

\begin{abstract}
For every positive integer $a$ which is coprime with $p$, $p$ is an odd prime, we denote
by $\overline{a}$ the unique integer satisfying $1\leq  \overline{a}\leq   p$ and $a\overline{a}\equiv 1(\mathrm{mod}~p)$. Put
$$L(p)=\{a\in Z^+:(a,p)=1,2\nmid a+\overline{a}\}.$$
The elements of $L(p)$ are called D.H. Lehmer numbers. The main purpose of this paper is to prove that for any fixed odd prime $p$, every sufficiently large number unless it is congruent to 15 or 16$(\mathrm{mod}~{16})$ is representable as the sum of 14 fourth powers of D.H. Lehmer numbers. Furthermore, every sufficiently large number is representable as the sum of 16 fourth powers of D.H. Lehmer numbers.

\end{abstract}

\footnote {2020 Mathematics Subject Classification{: Primary 11P05; Secondary 11P55.}}

\maketitle

\section{Introduction}
Let $m \geq 3$ be an odd integer. The classical problem of the D.H. Lehmer studies the properties of
$$r(m) = \underset{\substack{ab \equiv 1(\mathrm{mod}~m)\\ 2 \nmid a+b}}{\sum_{\substack{a=1 }}^{m}{}^{\prime} \sum_{\substack{b=1 }}^{m}{}^{\prime}}\: 1,$$
where $\sum_{a=1}^{m}{}^{\prime}$ means that the sum runs over the integers coprime with $m$. \\

From 1993 to 1994, Zhang Wenpeng\textsuperscript{\cite{5,6}} gave an asymptotic formula for $r(m)$, which reads as follows:
$$
r(m) = \frac{1}{2}\phi(m) + O(m^{\frac{1}{2}} d^2(m) \log^2 m).$$
In 1993, Zheng Zhiyong \textsuperscript{\cite{7}} considered short intervals, in fact, he proved that
$$\sum_{\substack{a \leq   N \\ a \in L(m)}} 1 = \frac{1}{2}N\phi(m)m^{-1} + O(m^{\frac{1}{2}} d^2(m) \log^2 m),$$
where $L(m) = \{a \in \mathbb{Z}^+, a \leq   m, (a,m)=1, a+\bar{a} \equiv 1 \pmod 2\}$ and $\bar{a}$ is the unique integer satisfying $\bar{a}a \equiv 1 \pmod m$, $1 \leq   \bar{a} \leq   m$.\\

For $n \geq 2$ be a fixed positive integer, and let $m \geq 3$ and $c$ be two integers with $(n,m)=(c,m)=1$. Denote by
$$
r_n(\delta_1,\delta_2,c;m) = \underset{\substack{ ab \equiv c(\mathrm{mod}~m) \\ n \nmid a+b}}{\sum_{\substack{a \leq   \delta_1 m \\  }}{}^{\prime} \sum_{\substack{b \leq   \delta_2 m}}{}^{\prime}}\:1  \quad(0 < \delta_1, \delta_2 \leq   1),$$
where $\sum'_{a\leq  \delta_1 m}$ means $\sum_{\substack{a\leq   \delta_1 m\\(a,m)=1}}$. In 2008, Lu Yaming and Yi Yuan\textsuperscript{\cite{12}} proved that
$$r_n(\delta_1,\delta_2,c;m)=\bigg(1-\frac{1}{n}\bigg)\delta_1\delta_2\phi(m)+O(m^{\frac{1}{2}}d^6(m)\log^2 m).$$
And in the same year they\textsuperscript{\cite{2}} studied the application of  D.H. Lehmer numbers in the Three Prime Theorem to arrive at the following conclusion, every integer $N$ like above can be expressed as the sum of three D.H. Lehmer numbers. Denote by $R(N)$ the number of the representations of $N$ as the sum of three D.H. Lehmer numbers, then
$$R(N) = \frac{N^2}{2} \bigg(1 - \frac{1}{n}\bigg)^3 \frac{\phi^3(m)}{m^3} A(m, N) + O(N^2 m^{-\frac{1}{2}} d^9(m) \log^3 m),$$
where
$$
A(m, N) = \prod_{p \mid (m, N)} \bigg(1 - \frac{1}{(p-1)^2}\bigg) \prod_{p \mid m, \, p \nmid N} \bigg(1 + \frac{1}{(p-1)^3}\bigg).$$\\
So similarly, are there similar results on the Waring's problem. \\

Write a sufficiently large positive integer $N$ as the sum of $s$ $k$th powers, i.e.
$$N=n_1^{k}+n_2^{k}+n_3^{k}+ \cdots +n_s^{k},$$
where $n_1,n_2,\cdots,n_s$ are positive integers. Defined for $k\geq 2$ to be the least such that every sufficiently large natural number is the sum of at most $s$ $k$th powers of natural numbers. difficult. In fact the value of $G(k)$ is only known when $k=2$ or $4$, namely,
$$G\left(2\right)=4,G\left(4\right)=16.$$
In 1943, U. V. Linnik\textsuperscript{\cite{4}} had shown that $G(3)\leq  7$ and in 1951 G. L. Watson\textsuperscript{\cite{8}} had given an extremely elegant proof of this. When $k>3$ all the best estimates available at present for $G(k)$ have been obtained via the Hardy-Littlewood method. \\
In 1939, H. Davenport\textsuperscript{\cite{1}} had already proven that every sufficiently large number can be representable as the sum of 14 fourth powers, unless it is congruent to 15 or 16$(\mathrm{mod}~16)$. And further proved that every sufficiently large number is representable as the sum of 16 fourth powers. And in fact, he\textsuperscript{\cite{3}} also proved that $G(4)=16$.\\
The authors are very interested in these problems. Since the relevant methods for the case where $k=4$ are quite mature, we hope to obtain the case where sufficiently large positive integers are written as fourth powers of D.H. Lehmer numbers. We will use the properties of the D.H. Lehmer numbers and the Hardy-Littlewood method to prove that every sufficiently large positive integer that meet the conditions is representable as the sum of 14 fourth powers of D.H. Lehmer numbers and every sufficiently large positive integer is representable as the sum of 16 fourth powers of D.H. Lehmer numbers. No one has yet conducted research on this problem, at least the authors have not see any relevant references. Therefore, we have the following theorem.\\

\noindent\textbf{Theorem 1.} \textit{Every sufficiently large positive integer unless it is congruent to 15 or 16 $(\mathrm{mod}~16)$ is representable as the sum of 14 fourth powers of the numbers in $L(p)$, where $L(p) = \{a \in \mathbb{Z}^+:  (a,p)=1,\bar{a}a \equiv 1(\mathrm{mod}~p), 1 \leq   \bar{a} \leq   p, a+\bar{a} \equiv 1 \pmod 2\}$, $p$ is any fixed odd prime. }\\

Obviously, $2p-1$ itself is a Lehmer number in $L(p)$, therefore, for $N\equiv15,16(\mathrm{mod}~16)$, we have $N-2(2p-1)^4\not \equiv15,16(\mathrm{mod}~16)$, then by Theorem 1, there are 14 numbers in $L(p)$ such that $$N-2(2p-1)^4=n_1^4+n_2^4+\cdots+n_{14}^4,$$
$$N=n_1^4+n_2^4+\cdots+n_{14}^4+(2p-1)^4+(2p-1)^4.$$ That is, we have the following Theorem 2.\\
 
\noindent\textbf{Theorem 2.} \textit{Every sufficiently large number is representable as the sum of 16 fourth powers of the numbers in $L(p)$, where $p$ is any fixed odd prime.}\\

\noindent\textbf{Note.} For the Hardy-Littlewood method commonly used in the Waring's problem, we have narrowed down the range of integers where $N$ is written as the sum of 16 fourth powers of integers to the D.H. Lehmer numbers in $L(p)$ of the original range $[P,2P]$. This also validates $G(4)=16$ from another perspective. In addition, for other cases of Waring's problem, the splitting of $N$ may be limited to the D.H. Lehmer numbers, this is only an open problem.
\section{Prerequisite Knowledge}
Firstly, we introduce some definitions to prepare for the proof of the theorem.\\

Let $N$ be a large integer, not congruent to 15 or 16$(\mathrm{mod}~16)$, which is to be represented as the sum of 14 fourth powers of the D.H. Lehmer numbers. We denote $\delta$ as a small positive number, fixed throughout the paper, and let $p$ be an odd prime.

Let
$$P=\bigg[\bigg(\frac{N}{100}\bigg)^{\frac{1}{4}}\bigg], Q=P^{3+\delta}$$
and
$$\mu=\frac{243}{1567}.$$

Let $u_1 < u_2 < \cdots < u_U < P^{\mu+3}$ be the numbers less than $P^{\mu+3}$ representable as the sum of 4 fourth powers of the numbers in $L(p)$. And $u_i$ congruent to $f (\bmod\ 16)$, $i=1,2,\cdots,U$, where $f=0,1,2,3,$ or 4. Thus
$$N-2f\equiv1,2,3,4,5\mathrm{~or~}6\mathrm{~(mod~}16).$$
For any real $\alpha$, we define
$$T(\alpha)\:=\:\sum_{\substack {{x\sim P} \\ { x \in L(p)}}}\:e(\alpha x^4),$$$$U(\alpha)=\sum_{h=1}^Ue(\alpha u_h),$$
$$r_{14}(N)=\int_{0}^{1}T^6(\alpha)U^2(\alpha)e(-N\alpha)d\alpha.$$
Throughout the paper, we write $\sum_x{'}$ for $\sum_{x,(x,p)=1}$, and $a$, $q$ are subject to $a \leq   q$, $(a, q) = 1$. Let
$$S_{a,q} = \sum_{x=1}^{q} e_q(ax^4), \quad (e_q(A) = e\bigg(\frac{A}{q}\bigg)),$$
for any real $\beta$,
$$
v_0(\beta) = \sum_{n=P^4}^{(2P)^4} \frac{1}{4} n^{-\frac{3}{4}} e(\beta n),$$
and
$$T^*(\alpha, a, q) = q^{-1} S_{a,q} v_0\bigg(\alpha - \frac{a}{q}\bigg).$$
Define $\mathfrak{M}(q,a)$  as the interval $[a/q-1/qQ,a/q+1/qQ]$, then write $\mathfrak{M}$ for the union of all $\mathfrak{M}(q,a)$ with $1\leq    a\leq    q\leq    P^{1/2}$ and  $(a,q)=1.$ Obviously, $\mathfrak{M}(q,a)$ are disjoint. Define $\mathfrak{m}$ as the complement of $\mathfrak {M} $ in $[1/qQ,1+1/qQ] $.
\section{Some Lemmas}
To prove the validity of the Theorem, we also introduce the following lemmas:\\

\noindent{\textbf{Lemma 1.}}\textit{Let}
$$S(q,a,b)=\sum_{x=1}^{q}e((ax^{k}+bx)q^{-1}).$$
\textit{Suppose that $(q,a)=1$. Then}
$$
S(q,a,b)\ll q^{1/2+\varepsilon}(q,b).$$
\textbf{Proof.} See [3], Lemma 4.1.\\

\noindent{\textbf{Lemma 2.}}\textit{Suppose that $X < Y$, $F''$ exists and is continuous on $[X, Y]$ and $F'$ is monotonic on $[X, Y]$. Let $H_1$, $H_2$ denote integers such that $H_1 \leq   F'(\alpha) \leq   H_2$ for every $\alpha$ in $[X, Y]$. Then}
$$\sum_{X < x \leq   Y} e(F(x)) = \sum_{h = H_1}^{H_2} \int_X^Y e(F(\alpha) - \alpha h) d\alpha + O(\log(2 + H)),$$
\textit{where $H = \max(| H_1| , | H_2| )$.}\\
\textbf{Proof.} See [3], Lemma 4.2.
\\

\noindent{\textbf{Lemma 3.}}\textit{ Suppose that $(a,q)=1$, $\alpha=\frac{a}{q}+\beta$, and $n$ is large enough, define}
$$v(\beta)=\sum_{\substack{x\leq   n}}\frac{1}{k}x^{1/k-1}e(\beta x),\quad v
_1(\beta)=\int_{0}^{n^{1/k}}e(\beta \gamma^{k}) d\gamma,$$
$$f(\alpha)=\sum_{\substack{x\leq   X\\x\in L(p)}}e(\alpha x^k),\quad S(q,a)=\sum_{x\leq   q}e\bigg(\frac{a}{q}x^k\bigg),\quad S(q,a,b)=\sum_{x\leq   q}e\bigg(\frac{a}{q}x^k+\frac{b}{q}x\bigg),$$$$ V^{*}(\alpha,q,a)=q^{-1}S(q,a)v\bigg(\alpha-\frac{a}{q}\bigg),\quad  V_1^{*}(\alpha,q,a)=q^{-1}S(q,a)v_1\bigg(\alpha-\frac{a}{q}\bigg).$$
\textit{For $| \beta|  \leq   q^{-1}X^{1-k-\delta}$, there is a constant $C$ and function $H_k(\alpha)$, then}
$$f(\alpha)- \frac{p-1}{2p}V^{*}(\alpha,q,a)-H_k(\alpha)\ll q^{\frac{1}{2}+\epsilon}$$
\textbf{Proof.}Let $X = n^{\frac{1}{k}}$, then we have
\begin{align}
    \begin{split}
        f(\alpha) &= \sum_{\substack{x \leq   X\\ x \in L(p)}} e(\beta x^{k}) \sum_{\substack{m=1 \\ m \equiv x (\text{mod } q)}}^{q} e\bigg(\frac{a}{q}m^{k}\bigg)\\&= q^{-1} \sum_{-q/2 < b \leq   q/2} S(q,a,b) F(b),
    \end{split}
\end{align}
where
$$F(b) = \sum_{\substack{x \leq   X\\x\in L(p)}} e\bigg(\beta x^{k} - \frac{b}{q}x\bigg).$$
From the definitions of $F(b) $ and $L(p)$, we have
\begin{align*}
         F(b)&= \sum_{\substack{x \leq   X\\x\in L(p)}} e\bigg(\beta x^{k} - \frac{b}{q}x\bigg)=\underset{\substack{xx' \equiv 1(\mathrm{mod}~p)\\ 2 \nmid x+x'}}{\sum_{\substack{x=1 }}^{X}{}^{\prime} \sum_{\substack{x'=1 }}^{p}{}^{\prime}}\:  e\bigg(\beta x^{k} - \frac{b}{q}x\bigg)\\ &=\underset{\substack{xx' \equiv 1(\mathrm{mod}~p)\\}}{\sum_{\substack{x=1 }}^{X}{}^{\prime} \sum_{\substack{x'=1 }}^{p}{}^{\prime}}\:  e\bigg(\beta x^{k} - \frac{b}{q}x\bigg)-\underset{\substack{xx' \equiv 1(\mathrm{mod}~p)\\ 2 \mid x+x'}}{\sum_{\substack{x=1 }}^{X}{}^{\prime} \sum_{\substack{x'=1 }}^{p}{}^{\prime}}\:  e\bigg(\beta x^{k} - \frac{b}{q}x\bigg)=S_1(\alpha)-S_2(\alpha).
\end{align*}
Then
\begin{align}
    \begin{split}
        S_1(\alpha)&=\underset{\substack{xx' \equiv 1(\mathrm{mod}~p)\\}}{\sum_{\substack{x=1 }}^{X}{}^{\prime} \sum_{\substack{x'=1 }}^{p}{}^{\prime}}\:  e\bigg(\beta x^{k} - \frac{b}{q}x\bigg)=\frac{1}{p-1}\sum_{x\leq   X}{}^{\prime}\sum_{x'\leq   p}{}^{\prime}e\bigg(\beta x^{k} - \frac{b}{q}x\bigg)\sum_{\chi\mathrm{mod}p}\chi(xx')\\&=\frac{1}{p-1}\sum_{x\leq   X}{}^{\prime}\sum_{x'\leq   p}{}^{\prime}e\bigg(\beta x^{k} - \frac{b}{q}x\bigg)+\frac{1}{p-1}\sum_{x\leq   X}{}^{\prime}\sum_{x'\leq   p}{}^{\prime}e\bigg(\beta x^{k} - \frac{b}{q}x\bigg)\sum_{\substack{\chi\mathrm{mod}p\\\chi\neq\chi_0}}\chi(xx')\\&=\sum_{x\leq   X}{}^{\prime}e\bigg(\beta x^{k} - \frac{b}{q}x\bigg),
    \end{split}
\end{align}
and
\begin{align*}
    S_2(\alpha)&=\underset{\substack{xx' \equiv 1(\mathrm{mod}~p)\\ 2 \mid x+x'}}{\sum_{\substack{x=1 }}^{X}{}^{\prime} \sum_{\substack{x'=1 }}^{p}{}^{\prime}}\:  e\bigg(\beta x^{k} - \frac{b}{q}x\bigg)=\frac{1}{p-1}\underset{\substack{ 2 \mid x+x'}}{\sum_{x\leq   X}{}^{\prime}\sum_{x'\leq   p}{}^{\prime}}e\bigg(\beta x^{k} - \frac{b}{q}x\bigg)\sum_{\chi\mathrm{mod}p}\chi(xx')\\
    &=\frac{1}{p-1}\underset{\substack{ 2 \mid x+x'}}{\sum_{x\leq   X}{}^{\prime}\sum_{x'\leq   p}{}^{\prime}}e\bigg(\beta x^{k} - \frac{b}{q}x\bigg)+\frac{1}{p-1}\underset{\substack{ 2 \mid x+x'}}{\sum_{x\leq   X}{}^{\prime}\sum_{x'\leq   p}{}^{\prime}}e\bigg(\beta x^{k} - \frac{b}{q}x\bigg)\sum_{\substack{\chi\mathrm{mod}p\\\chi\neq \chi_0}}\chi(xx')\\
    &=S_{21}(\alpha)-S_{22}(\alpha).
\end{align*}
For $S_{21}(\alpha),$
\begin{align}
    \begin{split}
        S_{21}(\alpha)&=\frac{1}{p-1}\underset{\substack{ 2 \mid x+x'}}{\sum_{x\leq   X}{}^{\prime}\sum_{x'\leq   p}{}^{\prime}}e\bigg(\beta x^{k} - \frac{b}{q}x\bigg)\\
        &=\frac{1}{p-1}{\sum_{x\leq   X}{}^{\prime}e\bigg(\beta x^{k} - \frac{b}{q}x\bigg)\sum_{\substack{x'\leq   p\\x'\equiv -x(\mathrm{mod}2)\\ (x',p)=1}}}1\\
        &=\frac{1}{2}\sum_{x\leq   X}{}^{\prime}e\bigg(\beta x^{k} - \frac{b}{q}x\bigg).
    \end{split}
\end{align}
For $S_{22}(\alpha)$, $\chi\neq \chi_0,$
$$\chi(a)=\frac{1}{q}\sum_{k=1}^{q}G(k,\chi)e\bigg(-\frac{ak}{q}\bigg)=\frac{1}{q}\sum_{k=1}^{q-1}G(k,\chi)e\bigg(-\frac{ak}{q}\bigg),$$
then
\begin{align*}
    S_{22}(\alpha)=&\frac{1}{p-1}\underset{\substack{ 2 \mid x+x'}}{\sum_{x\leq   X}{}^{\prime}\sum_{x'\leq   p}{}^{\prime}}e\bigg(\beta x^{k} - \frac{b}{q}x\bigg)\sum_{\substack{\chi\mathrm{mod}p\\\chi\neq \chi_0}}\chi(xx')\\=&\frac{1}{2p-1}{\sum_{x\leq   X}{}^{\prime}\sum_{x'\leq   p}{}^{\prime}}e\bigg(\beta x^{k} - \frac{b}{q}x\bigg)\sum_{l=1}^{2}{e\bigg(\frac{x+x'}{2}l\bigg)}\sum_{\substack{\chi\mathrm{mod}p\\\chi\neq \chi_0}}\chi(xx')\\=&\frac{1}{2p-1}\sum_{l=1}^{2}\sum_{\substack{\chi\mathrm{mod}p\\\chi\neq \chi_0}}{\sum_{x\leq   X}}\chi(x)e\bigg(\beta x^{k} - \frac{b}{q}x\bigg){e\bigg(\frac{x}{2}l\bigg)}\sum_{x'\leq   p}e\bigg(\frac{x'}{2}l\bigg)\chi(x')\\=&\frac{1}{2p^2p-1}\sum_{l=1}^{2}\sum_{\substack{\chi\mathrm{mod}p\\\chi\neq \chi_0}}\bigg(\sum_{k_1=1}^{p-1}G(k_1,\chi){\sum_{x\leq   X}}e(\beta x^{k} - \frac{b}{q}x-\frac{k_1}{p}x){e\bigg(\frac{x}{2}l\bigg)}\bigg)\\&\times\bigg(\sum_{k_2=1}^{p-1}G(k_2,\chi)\sum_{x'\leq   p}e\bigg(-\frac{k_2}{q}x'\bigg)e\bigg(\frac{x'}{2}l\bigg)\bigg)\\=&\frac{1}{2p^2p-1}\sum_{l=1}^{2}\sum_{k_2=1}^{p-1}\sum_{k_1=1}^{p-1}\bigg({\sum_{x\leq   X}}e\bigg(\beta x^{k}-\frac{b}{q}x-\frac{k_1}{p}x\bigg){e\bigg(\frac{x}{2}l\bigg)}\bigg)\\&\times\bigg(\sum_{x'\leq   p}e\bigg(-\frac{k_2}{q}x'\bigg)e\bigg(\frac{x'}{2}l\bigg)\bigg)\sum_{\substack{\chi\mathrm{mod}p\\\chi\neq \chi_0}}G(k_1,\chi)G(k_2,\chi)\\=&\sum_{l=1}^{2}\sum_{k_1=1}^{p-1}C_{p}(k_1,l)\bigg({\sum_{x\leq   X}}e\bigg(\beta x^{k}-\frac{b}{q}x-\frac{k_1}{p}x\bigg){e\bigg(\frac{x}{2}l\bigg)}\bigg),
\end{align*}
where $C_{p}(k_1,l)$ is $$C_{p}(k_1,l)=\frac{1}{2p^2p-1}\sum_{k_2=1}^{p-1}(\sum_{x'\leq   p}e\bigg(-\frac{k_2}{q}x'\bigg)e(\frac{x'}{2}l))\sum_{\substack{\chi\mathrm{mod}p\\\chi\neq \chi_0}}G(k_1,\chi)G(k_2,\chi),$$
Hence
\begin{align}
    \begin{split}
        S_{22}(\alpha)=&\sum_{k_1=1}^{p-1}C_{p}(k_1,2){\sum_{x\leq   X}}e\bigg(\beta x^{k}-\frac{b}{q}x-\frac{k_1}{p}x\bigg)\\&+\sum_{k_1=1}^{p-1}C_{p}(k_1,1){\sum_{x\leq   X}}(-1)^xe\bigg(\beta x^{k}-\frac{b}{q}x-\frac{k_1}{p}x\bigg)\\=&\sum_{k_1=1}^{p-1}C_{p}(k_1,2){\sum_{x\leq   X}}e\bigg(\beta x^{k}-\frac{b}{q}x-\frac{k_1}{p}x\bigg)\\&+2\sum_{k_1=1}^{p-1}C_{p}(k_1,1)\underset{2\mid x}{\sum_{x\leq   X}}e\bigg(\beta x^{k}-\frac{b}{q}x-\frac{k_1}{p}x\bigg)\\&-\sum_{k_1=1}^{p-1}C_{p}(k_1,1){\sum_{x\leq   X}}e\bigg(\beta x^{k}-\frac{b}{q}x-\frac{k_1}{p}x\bigg).
    \end{split}
\end{align}
Then by (3.2), (3.3), (3.4), we have
\begin{align*}
    \begin{split}
        F(b)=&S_1(\alpha)-S_{21}(\alpha)-S_{22}(\alpha)\\=&\frac{1}{2}\sum_{x\leq   X}{}^{\prime}e\bigg(\beta x^{k} - \frac{b}{q}x\bigg)-\sum_{k_1=1}^{p-1}C_{p}(k_1,2){\sum_{x\leq   X}}e\bigg(\beta x^{k}-\frac{b}{q}x-\frac{k_1}{p}x\bigg)\\&-2\sum_{k_1=1}^{p-1}C_{p}(k_1,1)\underset{2\mid x}{\sum_{x\leq   X}}e\bigg(\beta x^{k}-\frac{b}{q}x-\frac{k_1}{p}x\bigg)\\&+\sum_{k_1=1}^{p-1}C_{p}(k_1,1){\sum_{x\leq   X}}e\bigg(\beta x^{k}-\frac{b}{q}x-\frac{k_1}{p}x\bigg)\\=&\frac{1}{2}\sum_{x\leq   X}e\bigg(\beta x^{k} - \frac{b}{q}x\bigg)-\frac{1}{2}\sum_{\substack{x\leq   X\\p\mid x}}e\bigg(\beta x^{k} - \frac{b}{q}x\bigg)\\&-\sum_{k_1=1}^{p-1}C_{p}(k_1,2){\sum_{x\leq   X}}e\bigg(\beta x^{k}-\frac{b}{q}x-\frac{k_1}{p}x\bigg)\\&-2\sum_{k_1=1}^{p-1}C_{p}(k_1,1)\underset{2\mid x}{\sum_{x\leq   X}}e\bigg(\beta x^{k}-\frac{b}{q}x-\frac{k_1}{p}x\bigg)\\&+\sum_{k_1=1}^{p-1}C_{p}(k_1,1){\sum_{x\leq   X}}e\bigg(\beta x^{k}-\frac{b}{q}x-\frac{k_1}{p}x\bigg).
    \end{split}
\end{align*}
By Lemma 2, the functions of $x$ in the exponential sums are naturally continuous on $[0,X]$, set $H=p^4,H_1=-H_2=-H$. Define $I(c)=\int_{0}^{X}e(\beta\gamma^k-\frac{c}{q}\gamma)d\gamma$, then
\begin{align}
    \begin{split}
        F(b)=&\frac{1}{2}\sum_{h=-H}^{H}I(b+hq)-\frac{1}{2p}\sum_{h=-H}^{H}I\bigg(b+\frac{hq}{p}\bigg)\\
        &-\sum_{k_1=1}^{p-1}(C_{p}(k_1,2)-C_{p}(k_1,1))\sum_{h=-H}^{H}I\bigg(b+hq+\frac{qk_1}{p}\bigg)\\
        &-2\sum_{k_1=1}^{p-1}C_{p}(k_1,1)\sum_{h=-H}^{H}I\bigg(b+\frac{qh}{2}+\frac{qk_1}{p}\bigg)+O(\mathrm{log}(2+H)).
    \end{split}
\end{align}\\
When $ b=0$, we note that the value within $I(c)$ in $F(b)$ is equal to $0$ if and only if $h=0$, and it only holds for $b+hq$ and $b+\frac{hq}{p}$. When $ b\neq0$, by (3.1), (3.5),
\begin{align*}
    \begin{split}
        f(\alpha)-\frac{p-1}{2p}V_1^*(\alpha)
        = &\frac{1}{2q} \underset{b+qh\neq0}{\sum_{-q/2 < b \leq   q/2}\sum_{h=-H}^{H}} S(q,a,b)I(b+qh)\\
        &-\frac{1}{2qp} \underset{b+\frac{hq}{p}\neq0}{\sum_{-q/2 < b \leq   q/2}\sum_{h=-H}^{H}} S(q,a,b)I\bigg(b+\frac{hq}{p}\bigg)\\
        &-\frac{1}{2qp} \underset{\substack{b+\frac{hq}{p}=0\\b\neq0}}{\sum_{-q/2 < b \leq   q/2}\sum_{h=-H}^{H}} S(q,a,b)I\bigg(b+\frac{hq}{p}\bigg)\\
        -\frac{1}{q}&\sum_{k_1=1}^{p-1}(C_{p}(k_1,2)-C_{p}(k_1,1))\underset{b+qh+\frac{qk_1}{p}\neq0}{\sum_{-q/2<b<q/2}\sum_{h=-H}^{H}}S(q,a,b)I\bigg(b+hq+\frac{qk_1}{p}\bigg)\\
        -\frac{1}{q}&\sum_{k_1=1}^{p-1}(C_{p}(k_1,2)-C_{p}(k_1,1))\underset{\substack{b+qh+\frac{qk_1}{p}=0\\b\neq0}}{\sum_{-q/2<b<q/2}\sum_{h=-H}^{H}}S(q,a,b)I\bigg(b+hq+\frac{qk_1}{p}\bigg)\\
        &-\frac{2}{q}\sum_{k_1=1}^{p-1}C_{p}(k_1,1)\underset{\substack{b+\frac{hq}{2}+\frac{qk_1}{p}\neq0}}{\sum_{-q/2<b<q/2}\sum_{h=-H}^{H}}S(q,a,b)I\bigg(b+\frac{qh}{2}+\frac{qk_1}{p}\bigg)\\
        &-\frac{2}{q}\sum_{k_1=1}^{p-1}C_{p}(k_1,1)\underset{\substack{b+\frac{hq}{2}+\frac{qk_1}{p}=0\\b\neq0}}{\sum_{-q/2<b<q/2}\sum_{h=-H}^{H}}S(q,a,b)I\bigg(b+\frac{qh}{2}+\frac{qk_1}{p}\bigg)\\
        &+O(\frac{1}{q}S(q,a,b)).
    \end{split}
\end{align*}
By Lemma 1, we have
\begin{align*}
    \begin{split}
       & f(\alpha)-\frac{p-1}{2p}V_1^*(\alpha)-H_k(\alpha)\\
        \ll &\frac{1}{2q} \underset{b+qh\neq0}{\sum_{-q/2 < b \leq   q/2}\sum_{h=-H}^{H}}q^{1/2+\epsilon}(q,b)\bigg| I(b+qh)\bigg| \\&+\frac{1}{2qp} \underset{b+\frac{hq}{p}\neq0}{\sum_{-q/2 < b \leq   q/2}\sum_{h=-H}^{H}}q^{1/2+\epsilon}(q,b)\bigg| I\bigg(b+\frac{hq}{p}\bigg)\bigg| \\&+\sum_{k_1=1}^{p-1}\bigg| C_{p}(k_1,2)+C_{p}(k_1,1)\bigg| \\&\times \underset{b+qh+\frac{qk_1}{p}\neq0}{\sum_{-q/2<b<q/2}\sum_{h=-H}^{H}}q^{1/2+\epsilon}(q,b)\bigg| I\bigg(b+hq+\frac{qk_1}{p}\bigg)\bigg| \\&+2\sum_{k_1=1}^{p-1}C_{p}(k_1,1)\underset{\substack{b+\frac{hq}{2}+\frac{qk_1}{p}\neq0}}{\sum_{-q/2<b<q/2}\sum_{h=-H}^{H}}q^{1/2+\epsilon}(q,b)\bigg| I\bigg(b+\frac{qh}{2}+\frac{qk_1}{p}\bigg)\bigg|\\&+O(q^{1/2+\epsilon}),
    \end{split}
\end{align*}
where
\begin{align*}
    H_k(\alpha)=&\frac{1}{q}\sum_{k_1=1}^{p-1}( C_{p}(k_1,2)-C_{p}(k_1,1)) \underset{\substack{b+qh+\frac{qk_1}{p}=0\\b\neq0}}{\sum_{-q/2<b<q/2}\sum_{h=-H}^{H}}S(q,a,b)v_1(\beta)\\&+\frac{2}{q}\sum_{k_1=1}^{p-1}C_{p}(k_1,1)\underset{\substack{b+\frac{hq}{2}+\frac{qk_1}{p}=0\\b\neq0}}{\sum_{-q/2<b<q/2}\sum_{h=-H}^{H}}S(q,a,b)v_1(\beta)\\&+\frac{1}{2qp} \underset{\substack{b+\frac{hq}{p}=0\\b\neq0}}{\sum_{-q/2 < b \leq   q/2}\sum_{h=-H}^{H}} S(q,a,b)v_1(\beta).
\end{align*}
Let
$$g(\gamma)=\beta\gamma^k-\frac{c}{q}\gamma.$$
For $c\neq0$, since $| \beta|  \leq   \frac{1}{qX^{k-1+\delta}}$, then $k\beta X^{k-1}\leq   \frac{k}{qX^{\delta}}$. As long as $X$ is large enough, we can ensure $\frac{k}{qX^{\delta}}\ll \frac{c}{2q}$, then
$$| g^{\prime}(\gamma)| =\bigg|(\beta\gamma^k-\frac{c}{q}\gamma)^{\prime}\bigg|=\bigg| k\beta\gamma^{k-1}-\frac{c}{q}\bigg| \gg \bigg| \frac{c}{2q}\bigg|,$$
furthermore
$$I(c)\ll\bigg| \frac{2q}{c}\bigg| .$$
Hence
\begin{align*}
    \begin{split}
       &f(\alpha)-\frac{p-1}{2p}V_1^*(\alpha)-H_k(\alpha)\\
       \ll &\frac{1}{2q} \underset{b+qh\neq0}{\sum_{-q/2 < b \leq   q/2}\sum_{h=-H}^{H}} q^{1/2+\epsilon}(q,b+qh)\bigg| \frac{b}{b+qh}\bigg| \\&+\frac{1}{2qp} \underset{b+\frac{hq}{p}\neq0}{\sum_{-q/2 < b \leq   q/2}\sum_{h=-H}^{H}} q^{1/2+\epsilon}(q,pb+qh)\bigg| \frac{2qp}{pb+qh}\bigg| \\&+\sum_{k_1=1}^{p-1}(\bigg| C_{p}(k_1,2)\bigg|+\bigg| C_{p}(k_1,1)\bigg|) \\&\times\underset{2bp+hqp+2qk_1\neq 0}{\sum_{-q/2<b<q/2}\sum_{h=-H}^{H}}q^{1/2+\epsilon}(q,bp+hqp+qk_1)\bigg| \frac{2qp}{bp+hqp+qk_1}\bigg| \\&+2\sum_{k_1=1}^{p-1}\bigg| C_{p}(k_1,1)\bigg| \underset{2bp+hqp+2qk_1\neq 0}{\sum_{-q/2<b<q/2}\sum_{h=-H}^{H}}q^{1/2+\epsilon}(q,2bp+hqp+2qk_1)\\&\times\bigg| 4qp/(2bp+hqp+2qk_1)\bigg| .
    \end{split}
\end{align*}
Through simulation,  $b+qh,pb+qh,bp+hqp+qk1,2bp+hqp+2qk_1$ in the above inequality, for a fixed $k_1$, different values of $b,h$ will result in a maximum of $p$ identical outcomes, since $q^{-1}\sum_{h=1}^q(q.b)\leq   d(q)$ and by (4.17) in [3] $V_1(\alpha,q,a)-V(\alpha,q,a)\ll 1$. Therefore
$$f(\alpha)-\frac{p-1}{2p}V_1^*(\alpha)-H_k(\alpha)\ll q^{1/2+\epsilon}\sum_{c=1}^{2p^5q}(q,c)\frac{q}{c}\ll q^{1/2+\epsilon}.$$
This proves Lemma 3.\\

\noindent{\textbf{Lemma 4.}} \textit{Let $c_1, c_2, \dots$ be any sequence of complex numbers and suppose that $F$ has a continuous derivative on $[0, X]$. Then}
$$\sum_{m \leq   X} c_m F(m) = F(X) \sum_{m \leq   X} c_m - \int_0^X F'(\gamma) \sum_{m \leq   \gamma} c_m \, d\gamma.$$
 \textbf{Proof.} See [3], Lemma 2.6.\\

\noindent{\textbf{Lemma 5.}} \textit{Suppose that $(a,q)=1$ and $\alpha\in \mathfrak{M}$, then }
$$f(\alpha)- \frac{p-1}{2p}V^{*}(\alpha,q,a)\leq  (1+| \beta| X)q.$$
 \textbf{Proof.} By the conclusion of [7], for $Y\geq 0$,
 $$\sum_{\substack{m\leq   Y\\m\in L(p)}}e\bigg(\frac{a}{q}m^{k}\bigg)=\sum_{r=1}^{q}e\bigg(\frac{a}{q}r^{k}\bigg)\sum_{\substack{m\leq   Y\\m\equiv r({\mathrm{mod}~p})\\m\in L(p)}}1=\frac{Y(p-1)}{2pq}S(q,a)+O(q),$$
and
$$\sum_{m\leq   Y^k}\frac{1}{k}m^{1/k-1}=\int_{1}^{Yk}\frac{1}{k}\alpha^{1/k-1}d\alpha+O(1)=Y+O(1).$$
Let
$$c_{m}=\begin{cases}
e(am/q)-\frac{p-1}{2pq}S(q,a)\frac{1}{k}m^{1/k-1},&\text{when $m$ is a $k$-th powers of D.H. Lehmer number;}\\
-\frac{p-1}{2pq}S(q,a)\frac{1}{k}m^{1/k-1},&\text{otherwise},\end{cases}$$
and take $Y=\gamma^{1/k}$. Then
$$\sum_{m\leq  \gamma}c_{m}\leq   (\gamma\geq0).$$
Therefore, by Lemma 4 with $G(\gamma)=e(\beta\gamma)$,
$$\sum_{m<X}c_{m}e(\beta m)=f(\alpha)- \frac{p-1}{2p}V^{*}(\alpha,q,a)\leq  (1+| \beta| X)q.$$
This proves Lemma 5.\\

\noindent{\textbf{Lemma 6 .}} \textit{Suppose that $P$ is a large positive integer, and that $u_1 < u_2 < \cdots < u_U < P^{\mu+3}$, where $0 < \mu \leq   \frac{1}{2}$. Suppose also that}
$$U > P^{3(1-\mu)-\epsilon}.$$
\textit{Then the number of solutions of}
\begin{equation}
    x^{4}+u_{h}=y^{4}+u_{i},
\end{equation}
\textit{subject to}
\begin{equation}
    P\leq   x\leq   2P,\quad P\leq   y\leq   2P,\quad x,y\in L(p)
\end{equation}
\textit{is}
$$O(P^{2}U^{2}P^{3\mu-4+2\epsilon}).$$
\noindent{\textbf{Proof.}} This result can be naturally obtained from Lemma 1 of [1]. In fact, the condition of this Lemma only adds $x,y\in L(p)$. Despite the additional condition, the result remains bounded by that of Lemma 1 of [1].\\

\noindent{\textbf{Lemma 7.}} \textit{Suppose that $s \geq  2$, and that $f$ is one of $0, 1, 2, \cdots, s$. For $n > n_0(\epsilon)$, there exist at least $n^{\gamma_s - \epsilon}$ numbers less than $n$ representable as the sum of $s$ fourth powers of the numbers in $L(p)$ and  $n$ congruent to $f \pmod{16}$, where}
\begin{equation}
    \gamma_2 = \frac{1}{2}, \quad \gamma_3 = \frac{19}{28}, \quad \gamma_4 = \frac{331}{412}, \cdots,  \gamma_s = \frac{3 + 13\gamma_{s-1}}{4(3 + \gamma_{s-1})}.
\end{equation}
\noindent{Proof.} For $s=2$, let $f = f_1 + f_2$, where $f_1 = 0$ or $1$ and $f_2 = 0$ or $1$. The number of pairs $x_1, x_2$ satisfying
$$x_1 , x_2 < (\tfrac{1}{2}n)^{\frac{1}{4}},\quad x_1 , x_2 \in L(p),\quad x_1 \equiv f_1 \pmod{2}, \quad x_2 \equiv f_2 \pmod{2}$$
is greater than $C n^{\frac{1}{4}}$, where $C$ is a positive absolute constant. Since the number of representations of an integer $m$ as $x_1^4 + x_2^4$ is $O(m^{\epsilon})$. Hence the number of numbers that meet the conditions when $s=2$ is at least$$\frac{C n^{\frac{1}{4}}}{O(n^{\epsilon})} > n^{\frac{1}{4}-2\epsilon}.$$
For $s\geq 3$, We assume that the assertion of the lemma is true for $s - 1$, with $\frac{1}{2} \leq   \gamma_{s-1} < 1$, and deduce that it is true for $s$, where $\gamma_s$ is given by the last formula of (3.8).
Let $f = f_1 + f_2$, where $f_1$ is one of $0, 1, \cdots, s - 1$, and $f_2$ is $0$ or $1$. Let
$$\mu = \frac{3(1 - \gamma_{s-1})}{3 + \gamma_{s-1}}.$$
Since $\frac{1}{2} \leq   \gamma_{s-1} < 1$, we have $0 < \mu < \frac{1}{2}$. Let $P = [\frac{1}{2}(\frac{1}{2}n)]^{\frac{1}{4}}$. Let $u_1 < u_2 < \cdots < u_\nu < P^{\mu+3}$ be the numbers less than $P^{\mu+3}$ representable as the sum of $s - 1$ fourth powers of numbers of $L(p)$ and congruent to $f_1 \pmod{16}$. By hypothesis
$$
U > P^{(\mu+3)\gamma_{s-1} - \epsilon}= P^{3(1 - \mu) - \epsilon}.$$
Let $r(m)$ denote the number of representations of $m$ as $x^4 + u_h$, where $P \leq   x \leq   2P$,  $x \equiv f_2 \pmod{2}$ and $x\in L(p)$. Then
$$\sum_m r(m) \geq \tfrac{p-1}{4p} P U.$$
Moreover,  the conditions of the equation represented by $\sum_mr^2(m)$ are more than Lemma 6, so the conclusion of Lemma 6 can be used to determine the number of solutions to the $\sum_mr^2(m)$.That is to say
$$\sum_mr^2(m)=O(P^2U^2P^{3\mu-4+2\epsilon}).$$
Hence the number of numbers less than $n$ representable as the sum of $s$ fourth powers and congruent to $f \pmod{16}$ is
$$(\sum_mr(m))^2\\\geq{\sum_{\substack{m\\r(m)>0}}}1\quad\geq\quad\frac{(\sum_m r(m))^2}{\sum_mr^2(m)}>P^{4-3\mu-3\epsilon}\\>n^{\gamma_s-\epsilon},$$
where
$$\gamma_{s}=\frac{1}{4}(4-3\mu)=\frac{3+13\gamma_{s-1}}{4(3+\gamma_{s-1})}.$$
This proves Lemma 7.\\

Lemma 6 and Lemma 7 have proven that for the selected $u_1 < u_2 < \cdots < u_U $ can be expressed by the sum of 4 fourth powers of the numbers in $L(p)$  in this article, the number $U$ is actually
$$U>P^{(\mu+3)\gamma_4-\epsilon},\quad\gamma_4=\frac{331}{412}.$$\\

\noindent{\textbf{Lemma 8 .}} \textit{For $q>P^{1-\delta}$, $T(\alpha)=O(P^{7/8+\delta+\epsilon})$\\}
\textbf{\noindent{Proof.}} From the definitions of $T(\alpha) $ and $L(p)$, we have
\begin{align*}
         T(\alpha)&= \sum_{\substack{x\sim P\\x\in L(p)}} e(\alpha x^4)=\underset{\substack{xx' \equiv 1(\mathrm{mod}~p)\\ 2 \nmid x+x'}}{\sum_{\substack{x\sim P}}{}^{\prime} \sum_{\substack{x'=1 }}^{p}{}^{\prime}}\:  e(\alpha x^4)\\ &=\underset{\substack{xx' \equiv 1(\mathrm{mod}~p)\\}}{\sum_{\substack{x\sim P}}{}^{\prime} \sum_{\substack{x'=1 }}^{p}{}^{\prime}}\:  e(\alpha x^4)-\underset{\substack{xx' \equiv 1(\mathrm{mod}~p)\\ 2 \mid x+x'}}{\sum_{\substack{x\sim P}}{}^{\prime} \sum_{\substack{x'=1 }}^{p}{}^{\prime}}\:  e(\alpha x^4)=T_1(\alpha)-T_2(\alpha).
\end{align*}
Then
\begin{align}
    \begin{split}
        T_1(\alpha)&=\underset{\substack{xx' \equiv 1(\mathrm{mod}~p)\\}}{\sum_{\substack{x\sim P}}{}^{\prime} \sum_{\substack{x'=1 }}^{p}{}^{\prime}}\:  e(\alpha x^4)=\frac{1}{p-1}\sum_{x\sim P}{}^{\prime}\sum_{x'\leq   p}{}^{\prime}e(\alpha x^4)\sum_{\chi\mathrm{mod}p}\chi(xx')\\&=\frac{1}{p-1}\sum_{x\sim P}{}^{\prime}\sum_{x'\leq   p}{}^{\prime}e(\alpha x^4)+\frac{1}{p-1}\sum_{x\sim P}{}^{\prime}\sum_{x'\leq   p}{}^{\prime}e(\alpha x^4)\sum_{\substack{\chi\mathrm{mod}p\\\chi\neq\chi_0}}\chi(xx')\\&=\sum_{x\sim P}{}^{\prime}e(\alpha x^4),
    \end{split}
\end{align}
and
\begin{align}
    \begin{split}
        T_2(\alpha)&=\underset{\substack{xx' \equiv 1(\mathrm{mod}~p)\\ 2 \mid x+x'}}{\sum_{\substack{x\sim P}}{}^{\prime} \sum_{\substack{x'=1 }}^{p}{}^{\prime}}\:  e(\alpha x^4)=\frac{1}{p-1}\underset{\substack{ 2 \mid x+x'}}{\sum_{x\sim P}{}^{\prime}\sum_{x'\leq   p}{}^{\prime}}e(\alpha x^4)\sum_{\chi\mathrm{mod}p}\chi(xx')\\
    &=\frac{1}{p-1}\underset{\substack{ 2 \mid x+x'}}{\sum_{x\sim P}{}^{\prime}\sum_{x'\leq   p}{}^{\prime}}e(\alpha x^4)+\frac{1}{p-1}\underset{\substack{ 2 \mid x+x'}}{\sum_{x\sim P}{}^{\prime}\sum_{x'\leq   p}{}^{\prime}}e(\alpha x^4)\sum_{\substack{\chi\mathrm{mod}p\\\chi\neq \chi_0}}\chi(xx')\\
    &=T_{21}(\alpha)-T_{22}(\alpha),
    \end{split}
\end{align}
where
\begin{align}
    \begin{split}
        T_{21}(\alpha)&=\frac{1}{p-1}\underset{\substack{ 2 \mid x+x'}}{\sum_{x\sim P}{}^{\prime}\sum_{x'\leq   p}{}^{\prime}}e(\alpha x^4)=\frac{1}{p-1}{\sum_{x\sim P}{}^{\prime}e(\alpha x^4)\sum_{\substack{x'\leq   p\\x'\equiv -x(\mathrm{mod}2)\\ (x',p)=1}}} 1\\&=\frac{1}{2}\sum_{x\sim P}{}^{\prime}e(\alpha x^4).
    \end{split}
\end{align}
Since for $\chi\neq \chi_0$
$$\chi(a)=\frac{1}{q}\sum_{k=1}^{q}G(k,\chi)e\bigg(-\frac{ak}{q}\bigg)=\frac{1}{q}\sum_{k=1}^{q-1}G(k,\chi)e\bigg(-\frac{ak}{q}\bigg),$$
and by Weyl's inequality and the result of the proof process of Lemma 5 in [2],
\begin{align}
    \begin{split}
        T_{22}(\alpha)=&\frac{1}{p-1}\underset{\substack{ 2 \mid x+x'}}{\sum_{x\sim P}{}^{\prime}\sum_{x'\leq   p}{}^{\prime}}e(\alpha x^4)\sum_{\substack{\chi\mathrm{mod}p\\\chi\neq \chi_0}}\chi(xx')\\=&\frac{1}{2(p-1)}{\sum_{x\sim P}{}^{\prime}\sum_{x'\leq   p}{}^{\prime}}e(\alpha x^4)\sum_{l=1}^{2}{e\bigg(\frac{x+x'}{2}l\bigg)}\sum_{\substack{\chi\mathrm{mod}p\\\chi\neq \chi_0}}\chi(xx')\\=&\frac{1}{2(p-1)}\sum_{l=1}^{2}\sum_{\substack{\chi\mathrm{mod}p\\\chi\neq \chi_0}}{\sum_{x\sim P}}\chi(x)e(\alpha x^4){e\bigg(\frac{x}{2}l\bigg)}\sum_{x'\leq   p}e\bigg(\frac{x'}{2}l\bigg)\chi(x')\\
        \end{split}
\end{align}
\begin{align*}
    \begin{split}=&\frac{1}{2p^2(p-1)}\sum_{l=1}^{2}\sum_{\substack{\chi\mathrm{mod}p\\\chi\neq \chi_0}}\bigg(\sum_{k_1=1}^{p-1}G(k_1,\chi){\sum_{x\sim P}}e\bigg(\alpha x^4-\frac{k_1}{p}x\bigg){e\bigg(\frac{x}{2}l\bigg)}\bigg)\\&\times\bigg(\sum_{k_2=1}^{p-1}G(k_2,\chi)\sum_{x'\leq   p}e\bigg(-\frac{k_2}{p}x'\bigg)e\bigg(\frac{x'}{2}l\bigg)\bigg)\\=&\frac{1}{2p^2(p-1)}\sum_{l=1}^{2}\sum_{k_2=1}^{p-1}\sum_{k_1=1}^{p-1}\bigg({\sum_{x\sim P}}e\bigg(\alpha x^{4}-\frac{k_1}{p}x\bigg){e\bigg(\frac{x}{2}l\bigg)}\bigg)\\&\times\bigg(\sum_{x'\leq   p}e\bigg(-\frac{k_2}{p}x'\bigg)e\bigg(\frac{x'}{2}l\bigg)\bigg)\sum_{\substack{\chi\mathrm{mod}p\\\chi\neq \chi_0}}G(k_1,\chi)G(k_2,\chi)\\\ll& p^{1/2}d(p^3)\mathrm{log}p P^{1+\epsilon}(q^{-1}+P^{-1}+qP^{-4})^{1/8}\ll P^{7/8+\delta/8+\epsilon}.
    \end{split}
\end{align*}

By (3.9), (3.10), (3.11), (3.12), we obtain
\begin{equation}
    T(\alpha)=T_{1}(\alpha)-T_{21}(\alpha)-T_{22}(\alpha)=\frac{1}{2}\sum_{x\sim P}{}^{\prime}e(\alpha x^4)+O(P^{7/8+\epsilon})
\end{equation}
where
$$\sum_{x\sim P}{}^{\prime}e(\alpha x^4)=\sum_{d\mid p}\mu(d)\sum_{\substack{x\sim P/d}}e(\alpha d^4x^4).$$
By Dirichlet' s theorem on diophantine approximation, there are coprime integers $b$, $r$ with $r\leq  8P^3d^{-3}$, $| d^{4}\alpha-b/r| \leq  \frac{1}{8}r^{-1}d^{3}P^{-3}$. By Weyl' s inequality, when $r>P/d$
$$\sum_{P/d<y\leq  2P/d}e(\alpha d^4y^4)\ll \bigg(\frac{P}{d}\bigg)^{7/8},$$
and when $r\leq   P/d$, by Theorem 4.1 of [3],
\begin{align}
      \sum_{P/d<y\leq  2P/d}\:e(\alpha d^4y^4)\ll r^{-\frac{1}{4}} \frac{P}{d}\bigg(1+\bigg(\frac{P}{d}\bigg)^4\bigg| \alpha d^4-\frac{b}{r}\bigg| \bigg)^{-\frac{1}{4}}+\bigg(\frac{P}{d}\bigg)^{\frac{1}{2}+\epsilon}.
\end{align}
Furthermore, when
$$r\leq  \bigg(\frac{P}{d}\bigg)^{1/2};\quad\bigg| \alpha d^{4}-\frac{b}{r}\bigg| \leq  \frac{1}{r}\bigg(\frac{d}{P}\bigg)^{\frac{7}{2}},$$
it can also achieve a result of (3.14).
Hence
$$\sum_{x\sim P}{}^{\prime}e(\alpha x^4)\ll \sum_{\substack{d\mid p}}P^{7/8+\epsilon}d^{-7/8}+\sum_{\substack{d\mid p\\d\in \mathcal{D}}}\mu (d)Pd^{-1}r^{-\frac{1}{4}}\biggl(1+\biggl(\frac{P}{d}\biggr)^{4}\bigg| \alpha d^{4}-\frac{b}{r}\bigg| \biggr)^{-\frac{1}{4}},$$
where $\mathcal D $ =\{$| \alpha  d^4-b/r| \leq  {r^{-1}(P/d)^{-7/2}},r\leq  (P/d)^{1/2}$\}. Compare the conditions of $q,a$ and $b,r$, we have
$$\bigg| \frac{b}{r}-\frac{ad^4}{q}\bigg|\leq   \frac{1}{r}\bigg(\frac{d}{P}\bigg)^{\frac{7}{2}}+\frac{d^4}{qP^2},$$
then
$$| bq-ad^4r| \leq   qd^{\frac{7}{2}}P^{-\frac{7}{2}}+rP^{-2}d^4\ll 1,$$
since $P$ is large enough. Therefore $bq=ad^4r$, then $r=q/(q,d^4)$, by the trivial bound $(q,d^4)\leq   (q,d)^4,$
\begin{align}
    \begin{split}
        \sum_{x\sim P}{}^{\prime}e(\alpha x^4)&\ll P^{\frac{7}{8}+\epsilon}+q^{-1/4}P\bigg(1+P^4\bigg| \alpha-\frac{a}{q}\bigg| \bigg)^{-1/4}\sum_{d\mid p}\mu(d)\frac{(q,d)}{d}\\&\ll P^{7/8+\epsilon}+P^{3/4+\delta/4}\ll P^{7/8+\epsilon}.
    \end{split}
\end{align}
by (3.13) and (3.15), we have
$$T(\alpha)=O(P^{7/8+\delta+\epsilon}).$$
This proves Lemma 8.\\

\noindent{\textbf{Lemma 9 .}}\textit{Define}
$$S(q,\textbf{a})=S(q,a_{1},\ldots,a_{k})=\sum_{x=1}^{q}e((a_{1}x+\ldots+a_{k}x^{k})/q).$$
\textit{Suppose that} $(q,a_{1},\ldots,a_{k})=1$. \textit{Then}$$S(q,\textbf{a})\ll q^{1-1/k+\epsilon}.$$
\noindent{\textbf{Proof.}} See [3], Theorem 7.1.\\

\noindent{\textbf{Lemma 10 .}}\textit{ If $\mid \beta\mid \leq   \frac{1}{2}$, then}
$$v(\beta)=O(\min(P,P^{1-k}\mid\beta\mid^{-1})).$$
\noindent{\textbf{Proof.}} See [1], Lemma 4.\\

\noindent{\textbf{Lemma 11 .}} \textit{If $\alpha=\frac{a}{q}+\beta$, where$\mid\beta\mid\leq  \frac{1}{2}$, then }
$$T^{*}(\alpha,a,q)=O(q^{-1/4}\min(P,P^{-3}\mid\beta\mid^{-1})).$$
\noindent{\textbf{Proof.}} See [1], Lemma 5.\\

\noindent{\textbf{Lemma 12 .}} \textit{Write $H(\alpha)$ as $H_4(\alpha)$ in Lemma 3. For $q\leq   P^{1-\delta}$,}
$$H(\alpha)\ll q\text{,}\quad H(\alpha)\ll q^{-1/4}\min(P,P^{-3}\mid\beta\mid^{-1}).$$
\noindent{\textbf{Proof.}} By Lemma 3 and Lemma 5, we have
$$H(\alpha)+O(q^{1/2})\leq   (1+| \beta| P)q\ll q,$$
hence
$$H(\alpha)\ll q.$$
Moreover, by Lemma 9 and Lemma 10, we have
$$H(\alpha)\ll q^{-1}| S(q,\textbf{a})|| v(\beta)+O(1)| \ll q^{-1/4}\min(P,P^{-3}\mid\beta\mid^{-1}).$$This proves Lemma 12.\\

\noindent{\textbf{Lemma 13 .}} \textit{If $\alpha=\frac{a}{q}+\beta$, then }
$$\sum_{q\leq   P^{\frac{1}{2}}}\sum_a\int_{\mathfrak{M}_{a,q}}\mid T^6(\alpha)-(\frac{p-1}{2p}T^{*}(\alpha,a,q)+H(\alpha))^6\mid\mid U(\alpha)\mid ^2d\alpha=O(U^2P^{2-\frac{1}{4}+\epsilon}).$$
\noindent{\textbf{Proof.}} By Lemma 3, Lemma 11 and Lemma 12,
\begin{align*}
    \begin{split}
        &\int_{\mathfrak{M}_{a,q}}\mid T^6(\alpha)-(\frac{p-1}{2p}T^{*}(\alpha,a,q)+H(\alpha))^6\mid d\alpha\\
    \ll & q^{-\frac{3}{4}+\epsilon}\int_0^\infty\min(P^5,P^{-15}\beta^{-5})d\beta+q^{4+1/2}Q^{-1}\\
    \ll & q^{-\frac{3}{4}+\epsilon}P+q^{4+1/2}Q^{-1},
    \end{split}
\end{align*}
hence
\begin{align*}
    \begin{split}
        &\sum_{q\leq   P^{\frac{1}{2}}}\sum_a\int_{\mathfrak{M}_{a,q}}\mid T^6(\alpha)-(\frac{p-1}{2p}T^{*}(\alpha,a,q)+H(\alpha))^6\mid \mid U(\alpha)\mid^2d\alpha\\
        \ll & U^2\sum_{q\leq   P^{\frac{1}{2}}}(q^{\frac{1}{4}+\epsilon}P+q^{5+1/2}Q^{-1})\\
        \ll & U^2P^{1+\frac{1}{2}(\frac{5}{4}+\epsilon)}.
    \end{split}
\end{align*}

This proves Lemma 13.\\

\noindent{\textbf{Lemma 14 .}} \textit{For $P^{1/2}<q\leq   P^{1-\delta}$, $T(\alpha)=O(P^{7/8}).$\\}

\noindent{\textbf{Proof.}} By Lemma 3, Lemma 11 and Lemma 12,
\begin{align*}
    \begin{split}
        T(\alpha)&\leq   | T(\alpha)-T^{*}(\alpha,a,q)-H(\alpha)| +|T^{*}(\alpha,a,q)+H(\alpha)|\\
        &\ll q^{1/2+\epsilon}+q^{-1/4}\min(P,P^{-3}\mid\beta\mid^{-1})\\
        &\ll q^{1/2+\epsilon}+q^{-1/4}P\ll P^{7/8}.
    \end{split}
\end{align*}
This proves Lemma 14.\\

\noindent{\textbf{Lemma 15 .}} \textit{For $0\leq   l\leq   5$ is an integer, }
$$\sum_{q\leq   P^{\frac{1}{2}}}\sum_{a}\int_{\mathfrak{M}_{a,q}}\mid H(\alpha)^{6-l}(\frac{p-1}{2p}T^{*}(\alpha,q,a))^l\mid\mid U(\alpha)\mid ^{2}d\alpha=O(U^2P^{1+7/8+\epsilon}).$$
\noindent{\textbf{Proof.}} By Lemma 11 and Lemma 12,
\begin{align*}
    &\sum_{q\leq P^{\frac{1}{2}}}\sum_{a}\int_{\mathfrak{M}_{a,q}}\mid H(\alpha)^{6-l}T^{*l}(\alpha,q,a)\mid\mid U(\alpha)\mid ^{2}d\alpha\\
    \ll & U^2\sum_{q\leq   P^{\frac{1}{2}}}\sum_{a}q^{6-l-l/4}\int_0^\infty\min(P^l,P^{-3l}\beta^{-l})d\beta\\
    \ll & U^2P^{(8-l)/2+l-4}\\
    \ll & U^2P^{1+7/8+\epsilon}.
\end{align*}
This proves Lemma 15.\\

\noindent{\textbf{Lemma 16 .}} \textit{If $\overline{m}_{a,q}$ denotes the part of the interval $\frac{1}{Q}<\alpha<1+\frac{1}{Q}$ not belonging to $m_{a,q}$, then}
$$\sum_{q \leq   P^{\frac{1}{2}}} \sum_a \int_{\overline{m}_{a,q}} \mid T^*(\alpha, a, q)\mid ^6 \mid U(\alpha)\mid ^2 d\alpha = O(U^2).$$\\
\noindent{\textbf{Proof.}} See [1], Lemma 13.\\

\noindent{\textbf{Lemma 17 .}} \textit{There exists a $C_0>0$, we have}
$$\sum_{q\leq   P^{\frac{1}{2}}}\sum_a\int_0^1\bigg(\frac{p-1}{2p}T^{*}(\alpha,a,q)\bigg)^6U^2(\alpha)e(-N\alpha)d\alpha> C_0P^2U^2.$$
\noindent{\textbf{Proof.}} See the last part of [1]. In fact, the $T^*(\alpha,q,a)$ in this article is the same as $T^*(\alpha,q,a)$ in [1], so the result is only multiplied by an additional coefficient $(\frac{p-1}{2p})^6$, then we proves Lemma 17.\\

\noindent{\textbf{Lemma 18 .}} \textit{If  $\alpha=\frac{a}{q}+\beta$, where$\mid\beta\mid\leq  \frac{1}{2}$, then }
$$\int_{\mathfrak{m}}\mid T(\alpha)\mid^6\mid U(\alpha)\mid^2d\alpha=O(U^2P^{2-\frac{1}{2}+3\mu+5\delta}).$$
\noindent{\textbf{Proof.}} By Lemma 8 and Lemma 14, we have
$$T^4(\alpha)=O(P^{4-\frac{1}{2}+4\delta+4\epsilon}).$$
Moreover, for$$\int\mid T(\alpha)U(\alpha)\mid^2d\alpha,$$
taken over any interval of length 1, is precisely the number of solutions of (3.6) subject to (3.7), such that by Lemma 6,
$$\int\mid T(\alpha)U(\alpha)\mid^2d\alpha=O(P^{2}U^{2}P^{3\mu-4+2\epsilon}).$$
Hence $$\int_{\mathfrak{m}}\mid T(\alpha)\mid^6\mid U(\alpha)\mid^2d\alpha\ll T^4(\alpha)\int_0^1\mid T(\alpha)U(\alpha)\mid^2d\alpha
\ll U^2P^{2-\frac{1}{2}+3\mu+5\delta}.$$
This proves Lemma 18.\\

\section{Proof of the Theorem}
In this section, we will prove Theorem 1.\\
For
\begin{align}
    \begin{split}
        r_{14}(N)=&\int_{0}^{1} T^6(\alpha)U^2(\alpha)e(-N\alpha)d\alpha\\=&\sum_{q\leq   P^{\frac{1}{2}}}\sum_{a}\int_{0}^{1}\bigg(\frac{p-1}{2p}T^{*}(\alpha,a,q)\bigg)^6U^2(\alpha)e(-N\alpha)d\alpha\\&-\sum_{q\leq   P^{\frac{1}{2}}}\sum_a\int_{\overline{\mathfrak{M}}a,q}\bigg(\frac{p-1}{2p}T^{*}(\alpha,a,q)\bigg)^6U^2(\alpha)e(-N\alpha)d\alpha\\&+\sum_{l=0}^{5}\sum_{q\leq   P^{\frac{1}{2}}}\sum_{a}\int_{\mathfrak{M}_{a,q}} H(\alpha)^{6-l}T^{*l}(\alpha,q,a) U^{2}(\alpha)d\alpha\\&+\sum_{q\leq   P^{\frac{1}{2}}}\sum_a\int_{\mathfrak{M}_{a,q}}(T^6(\alpha)-(\frac{p-1}{2p}T^{*}(\alpha,a,q)+H(\alpha))^6)U^2(\alpha)e(-N\alpha)d\alpha\\&+\int_{\mathfrak{m}}T^6(\alpha)U^2(\alpha)e(-N\alpha)d\alpha.
    \end{split}
\end{align}
By Lemma 13, Lemma 15 and Lemma 16, we have the second, third and fourth sums of (4.1) is $O(P^{2-1/8+\epsilon}U^2)$, and by Lemma 16, the final integral is $O(U^2P^{2-\frac{1}{2}+3\mu+5\delta})$. Furthermore, the last three items of (4.1) is $O(U^2P^{2-\frac{1}{2}+3\mu+5\delta}).$
By Lemma 17,
$$\sum_{q\leq   P^{\frac{1}{3}}}\sum_{a}\int_{0}^{1}(\frac{p-1}{2p})^6T^{*^{6}}(\alpha,a,q)U^2(\alpha)e(-N\alpha)d\alpha > C_0 U^2P^2.$$
Since $\mu=\frac{243}{1567}$, therefore $2-\frac{1}{2}+3\mu+5\delta<2$.Hence
$$r_{14}(N)\gg 1.$$
This proves Theorem 1.\\

\hfill
$\square$

\vskip 3mm

\vskip 8mm

\noindent{\bf Acknowledgement}.
This work is supported in part by Shaanxi Fundamental Science Research Project for Mathematics and Physics
(Grant No. 23JSY033) and Natural Science Basic Research Project of Shaanxi Province (2021JM-044).

\end{document}